\documentclass[11pt]{amsart}

\usepackage[english]{babel}
\usepackage[T1]{fontenc}
\usepackage[utf8]{inputenc}

\usepackage{mathpazo}
\usepackage{graphicx}
\usepackage{microtype}
\usepackage{enumerate}
\usepackage{fontawesome}

\usepackage{csquotes}
\usepackage[dvipsnames]{xcolor}
\usepackage{url}

\usepackage{breakurl}
\usepackage[breaklinks,bookmarks, colorlinks, citecolor=MidnightBlue!75!black, urlcolor=MidnightBlue!75!black, linkcolor=MidnightBlue!75!black]{hyperref}

\calclayout

\AtBeginDocument{%
   \def\MR#1{}
}

\usepackage{scalerel}

\usepackage{amsmath, amssymb, amsthm, amscd}
\usepackage{commath, mathtools, mathrsfs}
\usepackage{thmtools}
\usepackage{tikz, tikz-cd}
\usetikzlibrary{matrix, arrows, decorations.pathmorphing}

\theoremstyle{plain}

\newtheorem{theorem}{Theorem}[section]
\newtheorem{lemma}[theorem]{Lemma}
\newtheorem{proposition}[theorem]{Proposition}
\newtheorem{corollary}[theorem]{Corollary}

\newtheorem*{theorem*}{Theorem}
\newtheorem*{lemma*}{Lemma}
\newtheorem*{proposition*}{Proposition}
\newtheorem*{corollary*}{Corollary}

\declaretheorem[numbered=no,name=Theorem A]{TA}
\declaretheorem[numbered=no,name=Theorem B]{TB}
\declaretheorem[numbered=no,name=Theorem C]{TC}

\theoremstyle{definition}
\newtheorem{remark}[theorem]{Remark}
\newtheorem{definition}[theorem]{Definition}

\newtheorem*{conjecture*}{Conjecture}
\newtheorem*{remark*}{Remark}
\newtheorem*{definition*}{Definition}
\newtheorem*{observation*}{Observation}

\newcommand{\on}{\operatorname}

\newcommand{\A}{\mathbb{A}}

\newcommand{\h}{\on{h}}

\newcommand{\mc}{\mathcal}

\newcommand{\Pic}{\on{Pic}}

\newcommand{\Q}{\mathbb{Q}}

\newcommand{\s}{\subseteq}

\newcommand{\spec}{\on{Spec}}

\renewcommand{\sup}{\scaleto{\triangle}{0.5em}}

\newcommand{\xar}{\xrightarrow}

\renewcommand{\epsilon}{\varepsilon}
\renewcommand{\H}{\on{H}}
\renewcommand{\hom}{\on{Hom}}

\renewcommand{\O}{\mc{O}}
\renewcommand{\P}{\mathbb{P}}
\renewcommand{\phi}{\varphi}

\usepackage[giveninits=true,style=alphabetic,backend=biber,maxnames=99,maxalphanames=5,giveninits=true]{biblatex}
\bibliography{khit}

\usepackage{todonotes}

\title{A higher dimensional Hilbert irreducibility theorem}
\date{}
\subjclass[2010]{}

\author{Giulio Bresciani}
\address{Freie Universit\"at Berlin, Arnimallee 3, 14195, Berlin, Germany}
\email{gbresciani@math.fu-berlin.de}
\thanks{The author is supported by the DFG Priority Program "Homotopy Theory and Algebraic Geometry" SPP 1786}

\begin{document}

\begin{abstract}
	Assuming the weak Bombieri-Lang conjecture, we prove that a generalization of Hilbert's irreducibility theorem holds for families of geometrically mordellic varieties (for instance, families of hyperbolic curves). As an application we prove that, assuming Bombieri-Lang, there are no polynomial bijections $\Q\times \Q\to\Q$. 
\end{abstract}

\maketitle
\tableofcontents

\section{Introduction}

Serre reformulated Hilbert's irreducibility theorem as follows \cite[Chapter 9]{ser97}.

\begin{theorem*}[Hilbert's irreducibility, Serre's form]
	Let $k$ be finitely generated over $\Q$, and let $f:X\to\P^{1}$ be a morphism with $X$ a scheme of finite type over $k$. Suppose that the generic fiber is finite, and that there are no generic sections $\spec k(\P^{1})\to X$. Then $X(k)\to\P^{1}(k)$ is not surjective.
\end{theorem*}

Recall that the weak Bombieri-Lang conjecture states that, if $X$ is a positive dimensional variety of general type over a field $k$ finitely generated over $\Q$, then $X(k)$ is not dense in $X$.

A variety $X$ over a field $k$ is \emph{geometrically mordellic}, or \emph{GeM}, if every subvariety of $X_{\bar{k}}$ is of general type. This generalizes to defining a scheme $X$ as geometrically mordellic, or GeM, if it is of finite type over $k$ and every subvariety of $X_{\bar{k}}$ is of general type. If the weak Bombieri-Lang conjecture holds and $k$ is a field finitely generated over $\Q$, then the set of rational points of a GeM scheme over $k$ is finite, since its Zariski closure cannot have positive dimension. 

Assuming Bombieri-Lang, we prove that Hilbert's irreducibility theorem generalizes to morphisms whose generic fiber is GeM.

\begin{TA}
	Let $k$ be finitely generated over $\Q$, and let $f:X\to\P^{1}$ be a morphism with $X$ a scheme of finite type over $k$. Suppose that the generic fiber is GeM, and that there are no generic sections $\spec k(\P^{1})\to X$.
	
	Assume either that the weak Bombieri-Lang conjecture holds in every dimension, or that it holds up to dimension equal to $\dim X$ and that there exists an $N$ such that $|X_{v}(k)|\le N$ for every rational point $v\in \P^{1}(k)$. Then $X(k)\to\P^{1}(k)$ is not surjective.
\end{TA}

There is a version of Hilbert's irreducibility theorem over non-rational curves, and the same is true for the higher dimensional generalization.

\begin{TB}
	Assume that the weak Bombieri-Lang conjecture holds in every dimension. Let $k$ be finitely generated over $\Q$, and let $f:X\to C$ be a morphism with $X$ any scheme of finite type over $k$ and $C$ a geometrically connected curve. Assume that the generic fiber is GeM, and that there are no generic sections $\spec k(C)\to X$. Then $X(h)\to C(h)$ is not surjective for some finite extension $h/k$.
\end{TB}

As an application of Theorem A we give an answer to a long-standing Mathoverflow question \cite{zh} which asks whether there exists a polynomial bijection $\Q\times\Q\to\Q$, conditional on the weak Bombieri-Lang conjecture.

\begin{TC}
	Assume that the weak Bombieri-Lang conjecture for surfaces holds, and let $k$ be a field finitely generated over $\Q$. There are no polynomial bijections $k\times k\to k$.
\end{TC}

We remark that B. Poonen has proved that, assuming the weak Bombieri-Lang conjecture for surfaces, there are polynomials giving \emph{injective} maps $\Q \times \Q \to \Q$, see \cite{poo10}.

In 2019, T. Tao suggested on his blog \cite{tao19} a strategy to try to solve the problem of polynomial bijections $\Q\times\Q\to\Q$ conditional on Bombieri-Lang, let us summarize it. Given a morphism $\A^{2}\to\A^{1}$ and a cover $c:\A^{1}\dashrightarrow \A^{1}$, denote by $P_{c}$ the pullback of $\A^{2}$. If $P_{c}$ is of general type, by Bombieri-Lang $P_{c}(\Q)$ is not dense in $P_{c}$ and hence by Hilbert irreducibility a generic section $\A^{1}\dashrightarrow P_{c}$ exists. If $P_{c}$ is of general type for "many" covers $c$, one might expect this to force the existence a generic section $\A^{1}\dashrightarrow\A^{2}$, which would be in contradiction with the bijectivity of $\A^{2}(\Q)\to\A^{1}(\Q)$.

The strategy had some gaps, though. There were no results showing that the pullback $P_{c}$ is of general type for "many" covers $c$, and it was not clear how this would force a generic section of $\A^{2}\to\A^{1}$. Tao started a so-called "polymath project" in order to crowdsource a formalization. The project was active for roughly one week in the comments section of the blog but didn't reach a conclusion. Partial progress was made, we cite the two most important contributions. W. Sawin showed that $\A^{2}(\Q)\to\A^{1}(\Q)$ can't be bijective if the generic fiber has genus $0$ or $1$. H. Pasten showed that, for some morphisms $\A^{2}\to\A^{1}$ with generic fiber of genus at least $2$, the base change of $\A^{2}$ along the cover $z^{2}-b:\A^{1}\to\A^{1}$ is of general type for a generic $b$.

Theorem A is far more general than Theorem C, but it is possible to extract from the proof of the former the minimal arguments needed in order to prove the latter. These minimal arguments are a formalization of the ideas described above, hence as far as Theorem C is concerned we have essentially filled in the gaps in Tao's strategy.



\subsection*{Acknowledgements}

I would like to thank Hélène Esnault for reading an earlier draft of the paper and giving me a lot of valuable feedback, and Daniel Loughran for bringing to my attention the problem of polynomial bijections $\Q\times\Q\to\Q$.

\subsection*{Conventions}

A variety over $k$ is a geometrically integral scheme of finite type over $k$. A smooth, projective variety is of general type if its Kodaira dimension is equal to its dimension: in particular, a point is a variety of general type. 

We say that a variety is of general type if it is birational to a smooth, projective variety of general type. More generally, we define the Kodaira dimension of any variety $X$ as the Kodaira dimension of any smooth projective variety birational to $X$.

Curves are assumed to be smooth, projective and geometrically connected. Given a variety $X$ (resp. a scheme of finite type $X$) and $C$ a curve, a morphism $X\to C$ is a family of varieties of general type (resp. of GeM schemes) if the generic fiber is a variety of general type (resp. a GeM scheme). Given a morphism $f:X\to C$, a generic section of $f$ is a morphism $s:\spec k(C)\to X$ (equivalently, a rational map $s:C\dashrightarrow X$) such that $f\circ s$ is the natural morphism $\spec k(C)\to C$ (equivalently, the identity $C\dashrightarrow C$).

\section{Pulling families to maximal Kodaira dimension}

This section is of purely geometric nature, thus we may assume that $k$ is algebraically closed of characteristic $0$ for simplicity. The results then descend to non-algebraically closed fields with standard arguments.

Given a family $f:X\to\P^{1}$ of varieties of general type and $c:\P^{1}\to \P^{1}$ a finite covering, let $f_{c}:X_{c}\to \P^{1}$ be the fiber product and, by abuse of notation, $c:X_{c}\to X$ the base change of $c$. The goal of this section is to obtain sufficient conditions on $c$ such that $X_{c}$ is of general type. This goal will be reached in \autoref{pull}, which contains all the geometry we'll need for arithmetic applications.

Let us say that $X\to\P^{1}$ is birationally trivial if there exists a birational morphism $X\dashrightarrow F\times \P^{1}$ which commutes with the projection to $\P^{1}$. If $f$ is birationally trivial, then clearly our goal is unreachable, since $X_{c}$ will have Kodaira dimension $-\infty$ no matter which cover $c:\P^{1}\to\P^{1}$ we choose. We will show that this is in fact the only exception.

Assume that $X$ is smooth and projective (we can always reduce to this case), then the relative dualizing sheaf $\omega_{f}$ exists \cite[Corollary 24]{kle80}.  First, we show that for \emph{every} non-birationally trivial family there exists an integer $m$ such that $f_{*}\omega_{f}^{m}$ has \emph{some} positivity \ref{pos}. Second, we show that if $f_{*}\omega_{f}^{m}$ has \emph{enough} positivity, then $X$ is of general type \ref{posgen}. We then pass from "some" to "enough" positivity by base changing along a cover $c:\P^{1}\to\P^{1}$.

\subsection{Positivity of $f_{*}\omega_{f}^{m}$}

There are two cases: either there exists some finite cover $c:C\to\P^{1}$ such that $X_{d}\to C$ is birationally trivial, or not. Let us say that $f:X\to\P^{1}$ is \emph{birationally isotrivial} in the first case, and non-birationally isotrivial in the second case.

The non-birationally isotrivial case has been extensively studied by Viehweg and Koll\'ar, we don't need to do any additional work.

\begin{proposition}[{Koll\'ar, Viehweg \cite[Theorem p.363]{kol87}}]\label{nisopos}
	Let $f:X\to \P^{1}$ be a non-birationally isotrivial family of varieties of general type, with $X$ smooth and projective. There exists an $m>0$ such that, in the decomposition of $f_{*}\omega_{f}^{m}$ in a direct sum of line bundles, each factor has positive degree.\qed
\end{proposition}

We are thus left with studying the positivity of $f_{*}\omega_{f}^{m}$ in the birationally isotrivial, non-birationally trivial case. We'll have to deal with various equivalent birational models of families, not always smooth, so let us first compare their relative pluricanonical sheaves.

\subsubsection{Morphisms of pluricanonical sheaves}

In this subsection, fix a base scheme $S$. If a morphism to $S$ is given, it is tacitly assumed to be flat, locally projective, finitely presentable, with Cohen-Macauley equidimensional fibers of dimension $n$. For such a morphism $f:X\to S$, the relative dualizing sheaf $\omega_{f}$ exists and is coherent, see \cite[Theorem 21]{kle80}. Recall that $\omega_{f}$ satisfies the functorial isomorphism
\[f_{*}\underline{\hom}_{X}(F,\omega_{f}\otimes_{X}f^{*}N)\simeq \underline{\hom}_{S}(R^{n}f_{*}F,N)\]
for every quasi-coherent sheaf $F$ on $X$ and every quasi-coherent sheaf $N$ on $S$. Write $\omega_{f}^{\otimes m}$ for the $m$-th tensor power, we may drop the superscript $\_^{\otimes}$ and just write $\omega_{f}^{m}$ if $\omega_{f}$ is a line bundle.

Every flat, projective map $f:X\to S$ of smooth varieties over $k$ satisfies the above, see \cite[Corollary 24]{kle80}, and in this case we can compute $\omega_{f}$ as $\omega_{X}\otimes f^{*}\omega_{S}^{-1}$, where $\omega_{X}$ and $\omega_{S}$ are the usual canonical bundles. Moreover, the relative dualizing sheaf behaves well under base change along morphisms $S'\to S$, see \cite[Proposition 9.iii]{kle80}. 

Given a morphism $g:Y\to X$ over $S$ and a quasi-coherent sheaf $F$ over $Y$, then $R^{n}f_{*}(g_{*}F)$ is the $E^{n,0}_{2}$ term of the Grothendieck spectral sequence $(R^{p}f\circ R^{q}g)(F)\Rightarrow R^{p+q}(f\circ g)(F)$, thus there is a natural morphism $R^{n}f_{*}(g_{*}F)\to R^{n}(fg)_{*}F$. This induces a natural map
\[\hom_{Y}(F,\omega_{fg})=\hom_{S}(R^{n}(fg)_{*}F,\O_{S})\to\hom_{S}(R^{n}f_{*}(g_{*}F),\O_{S})=\hom_{X}(g_{*}F,\omega_{f}).\]

\begin{definition}
	If $g:Y\to X$ is a morphism over $S$, define $g_{\sup,f}:g_{*}(\omega_{fg})\to\omega_{f}$ as the sheaf homomorphism induced by the identity of $\omega_{fg}$ via the homomorphism
	\[\hom_{Y}(\omega_{fg},\omega_{fg})\to\hom_{X}(g_{*}\omega_{fg},\omega_{f})\]
	given above for $F=\omega_{fg}$. With an abuse of notation, call $g_{\sup,f}$ the induced sheaf homomorphism $g_{*}(\omega_{fg}^{\otimes m})\to\omega_{f}^{\otimes m}$ for every $m\ge 0$. If there is no risk of confusion, we may drop the subscript $\__{f}$ and just write $g_{\sup}$.
\end{definition}

The following facts are straightforward, formal consequences of the definition of $g_{\sup}$, we omit proofs.

\begin{lemma}\label{supbc}
	Let $g:Y\to X$ be a morphism over $S$ and $s:S'\to S$ any morphism, $f':X'\to S'$, $g':Y'\to X'$ the pullbacks to $S'$. By abuse of notation, call $s$ the morphisms $Y'\to Y$, $X'\to X$, too. Then
	\[g'_{\sup}=g_{\sup}|_{X'}\in\hom_{X'}(g'_{*}\omega_{f'g'},\omega_{f'})=\hom_{X'}(s^{*}g_{*}\omega_{fg},s^{*}\omega_{f}).\]\qed
\end{lemma}

\begin{lemma}
	For every quasi-coherent sheaf $F$ on $Y$, the natural map
	\[\hom_{Y}(F,\omega_{fg})\to\hom_{X}(g_{*}F,\omega_{f})\]
	constructed above is given by 
	\[\phi\mapsto g_{\sup}\circ g_{*}\phi:g_{*}F\to g_{*}\omega_{fg}\to\omega_{f}.\]\qed
\end{lemma}

\begin{corollary}\label{supcomp}
	Let $h:Z\to Y$, $g:Y\to X$ be morphisms over $S$. Then, for every $m\ge0$,
	\[g_{\sup}\circ g_{*}h_{\sup}=(gh)_{\sup}:gh_{*}\omega_{fgh}^{\otimes m}\to g_{*}\omega_{fg}^{\otimes m}\to \omega_{f}^{\otimes m}.\]\qed
\end{corollary}

\begin{corollary}\label{supeq}
	Let $g:Y\to X$ be a morphism over $S$. Suppose that a group $H$ acts on $Y,X,S$ and $g,f$ are $H$-equivariant. Then $g_{*}\omega_{fg}^{\otimes m},\omega_{f}^{\otimes m}$ are $H$-equivariant sheaves and $g_{\sup}:g_{*}\omega_{fg}^{\otimes m}\to\omega_{f}^{\otimes m}$ is $H$-equivariant.\qed
\end{corollary}

\begin{lemma}\label{reld}
	Let $g:Y\to X$ be a morphism over $S$. Assume that $Y,X$ are smooth varieties over a field $k$, and that $g$ is birational. Then $g_{\sup}$ is an isomorphism.
	\begin{proof}
		We have $\omega_{f}=\omega_{X}\otimes f^{*}\omega_{S}^{-1}$ and $\omega_{fg}=\omega_{Y}\otimes (fg)^{*}\omega_{S}^{-1}$. Moreover, $\omega_{Y}=g^{*}\omega_{X}\otimes\O_{Y}(R)$ where $R$ is some effective divisor whose irreducible components are contracted by $g$, hence $\omega_{fg}=g^{*}\omega_{f}\otimes\O_{Y}(R)$. Since $g_{*}\O_{Y}(mR)\simeq\O_{X}$, we have a natural isomorphism $g_{*}(\omega_{fg}^{m})\simeq \omega_{f}^{m}$ by projection formula. This is easily checked to correspond to $g_{\sup}$, which is then an isomorphism as desired.
	\end{proof}
\end{lemma}

\subsubsection{Birationally isotrivial families}

Let $C$ be a smooth projective curve and $f:X\to C$ a birationally isotrivial family of varieties of general type, and let $F/k$ be a smooth projective variety such that the generic fiber of $f$ is birational to $F$. Let $H$ be the finite group of birational automorphisms of $F$. The scheme of fiberwise birational isomorphisms $\on{Bir}(X/C,F)\to C$ restricts to an $H$-torsor on some non-empty open subset $V$ of $C$. The action of $H$ on $\on{Bir}(X/C,F)|_{V}$ is transitive on the connected components, thus they are all birational.

\begin{definition}
	In the situation above, define $b:B_{f}\to C$ as the smooth completion of any connected component of $\on{Bir}(X/C,F)|_{V}$, and $G_{f}\s H$ as the subgroup of elements mapping $B_{f}$ to itself. Let us call $B_{f}\to C$ and $G_{f}$ the \emph{monodromy cover} and the \emph{monodromy group} of $f$ respectively.
\end{definition}

We have that $B_{f}\to C$ is a $G_{f}$-Galois covering characterized by the following universal property: if $C'$ is a smooth projective curve with a finite morphism $c:C'\to C$, then $X_{c}\to C'$ is birationally trivial if and only if there exists a factorization $C'\to B_{f}\to C$.

\begin{proposition}\label{isopos}
	Let $f:X\to C$ be a birationally isotrivial family of varieties of general type, with $X$ smooth and projective. If $p\in B_{f}$ is a ramification point of the monodromy cover $b:B_{f}\to C$, then for some $m$ there exists an injective sheaf homomorphism $\O_{B_{f}}(p)\to f_{b*}\omega_{f_{b}}^{m}$.
	\begin{proof}
		The statement is equivalent to the existence of a non-trivial section of $\omega_{f_{b}}^{m}$ which vanishes on the fiber $X_{b,p}$. Let $F$ be as above, $G_{f}$ acts faithfully with birational maps on $F$. By equivariant resolution of singularities, we may assume that $G_{f}$ acts faithfully by isomorphisms on $F$. We have that $X$ is birational to $(F\times B_{f})/G_{f}$ where $G_{f}$ acts diagonally.
		
		By resolution of singularities, let $X'$ be a smooth projective variety with birational morphisms $X'\to X$, $X'\to (F\times B_{f})/G_{f}$: thanks to \autoref{reld} we may replace $X$ with $X'$ and assume we have a birational morphism $X\to (F\times B_{f})/G_{f}$. By equivariant resolution of singularities again, we may find a smooth projective variety $Y$ with an action of $G_{f}$, a birational morphism $g:Y\to X_{b}$ and a birational, $G_{f}$-equivariant morphism $y:Y\to F\times B_{f}$. Call $\pi:F\times B_{f}\to B_{f}$ the projection.
		
		\[\begin{tikzcd}
			Y	\dar[swap,near start,"y"]\rar["g"]		&	X_{b}\rar["b"]												&	X\dar						\\
			F\times B_{f}\ar[rr]\ar[dr, swap, "\pi"]	&																&	(F\times B_{f})/G_{f}\dar		\\
														&	B_{f}\rar["b"]\ar[uul, leftarrow, crossing over, swap, near end, "\pi y"]\ar[uu, crossing over, leftarrow, near end, swap, "f_{b}"]		&	C
		\end{tikzcd}\]
		
		Recall that we are trying to find a global section of $\omega_{f_{b}}^{m}$ that vanishes on $X_{b,p}$, where $p$ is a ramification point of $b$. Thanks to \autoref{reld}, we have that $\pi y_{*}\omega_{\pi y}^{m}\simeq \pi_{*}\omega_{\pi}^{m}\simeq\O_{B_{f}}\otimes\H^{0}(F,\omega_{F}^{m})$, thus $\H^0(Y,\omega_{\pi y}^{m})=\H^0(F,\omega_{F}^{m})=\H^{0}(Y_{p},\omega_{Y_{p}}^{m})$.
		
		The sheaf homomorphism $g_{\sup}=g_{\sup,f_{b}}:g_{*}\omega_{\pi y}^{m}\to\omega_{f_{b}}^{m}$ induces a linear map
		\[g_{\sup}(p):\H^0(Y_{p},\omega_{Y_{p}}^{m})=\H^0(Y,\omega_{\pi y}^{m})\xar{g_{\sup}}\H^0(X_{b},\omega_{f_{b}}^{m})\xar{\bullet|_{p}}\H^{0}(X_{b,p},\omega_{X_{b,p}}^{m})\]
		where the last map is the restriction to the fiber. Let $V\s B_{f}$ be the étale locus of $b:B_{f}\to C$. Since $X_{b}|_{V}$ is smooth, then $g_{\sup}$ restricts to an isomorphism on $X_{b}|_{V}$ thanks to \autoref{reld} and thus the map $\H^0(Y,\omega_{\pi y}^{m})\to\H^0(X_{b},\omega_{f_{b}}^{m})$ is injective. 
		
		We want to show that the restriction map $\H^0(X_{b},\omega_{f_{b}}^{m})\to\H^{0}(X_{b,p},\omega_{X_{b,p}}^{m})$ is not injective for some $m$, it is enough to show that $g_{\sup}(p)$ is not injective. Thanks to \autoref{supbc}, we have that $g_{\sup}(p)=g_{p,\sup}:\H^{0}(Y_{p},\omega_{Y_{p}}^{m})\to\H^{0}(X_{b,p},\omega_{X_{b,p}}^{m})$.
		
		Recall now that $G_{f}$ acts on $Y$. Let $G_{f,p}$ be the stabilizer of $p\in B_{f}$, it is a non-trivial group since $p$ is a ramification point. Thanks to \autoref{supeq}, the stabilizer $G_{f,p}$ acts naturally on $\H^{0}(Y_{p},\omega_{Y_{p}}^{m})$, $\H^{0}(F,\omega_{F}^{m})$, $\H^{0}(X_{b,p},\omega_{X_{b,p}}^{m})$, and the maps $y_{p,\sup}:\H^{0}(Y_{p},\omega_{Y_{p}}^{m})\simeq\H^{0}(F,\omega_{F}^{m})$, $g_{p,\sup}:\H^{0}(Y_{p},\omega_{Y_{p}}^{m})\to \H^{0}(X_{b,p},\omega_{X_{b,p}}^{m})$ are $G_{f,p}$-equivariant. Moreover, the action on $\H^{0}(X_{b,p},\omega_{X_{b,p}}^{m})$ is trivial since the action on $X_{b,p}$ is trivial. It follows that $g_{\sup}(p)$ is $G_{f,p}$-invariant, and hence to show that it is not injective for some $m$ it is enough to show that the action of $G_{f,p}$ on $\H^{0}(F,\omega_{F}^{m}))=\H^{0}(Y_{p},\omega_{Y_{p}}^{m})$ is not trivial for some $m$.
		
		Since $F$ is of general type, $F\dashrightarrow\P(\H^{0}(F,\omega_{F}^{m}))$ is generically injective for some $m$, fix it. Since the action of $G_{f,p}$ on $F$ is faithful, for every non-trivial $g\in G_{f,p}$ there exists a section $s\in\H^{0}(F,\omega_{F}^{m})$ and a point $v\in F$ such that $s(v)=0$ and $s(g(v))\neq 0$, in particular the action of $G_{f,p}$ on $\H^{0}(F,\omega_{F}^{m})$ is not trivial and we conclude.
	\end{proof}
\end{proposition}

\begin{corollary}\label{pos}
	Let $f:X\to \P^{1}$ be a non-birationally trivial family of varieties of general type, with $X$ smooth and projective. Then there exists an $m$ with an injective homomorphism $\O(1)\to f_{*}\omega_{f}^{m}$.
	\begin{proof}
		If $f$ is not birationally isotrivial, apply \autoref{nisopos}. Otherwise, $f$ is birationally isotrivial and not birationally trivial, thus the monodromy cover $b:B_{f}\to\P^{1}$ is not trivial. Since $\P^{1}$ has no non-trivial étale covers, we have that $B_{f}\to \P^{1}$ has at least one ramification point $p$. Let $m$ be the integer given by \autoref{isopos}, and write $f_{*}\omega_{f}^{m}=\bigoplus_{i}\O_{\P^{1}}(d_{i})$. Since $\O_{B_{f}}(p)\s f_{b*}\omega_{f_{b}}^{m}$ and $\omega_{f_{b}}=b^{*}\omega_{f}$, see \cite[Proposition 9.iii]{kle80}, there exists an $i$ with $d_{i}>0$.
	\end{proof}
\end{corollary}

\subsection{Pulling families to maximal Kodaira dimension}

Now that we have established a positivity result for $f_{*}\omega_{f}^{m}$ of any non-birationally trivial family $f:X\to\P^{1}$, let us use this to pull families to maximal Kodaira dimension.

\begin{proposition}\label{posgen}
	Let $f:X\to\P^{1}$ be a family of varieties of general type, with $X$ smooth and projective. Then $X$ is of general type if and only if there exists an injective homomorphism $\O_{\P^{1}}(1)\to f_{*}\omega_{X}^{m_{0}}$, or equivalently $\O_{\P^{1}}(2m_{0}+1)\to f_{*}\omega_{f}^{m_{0}}$, for some $m_{0}>0$.
	\begin{proof}
		By resolution of singularities, there exists a birational morphism $g:X'\to X$ with $X'$ smooth and projective such that the generic fiber of $X'\to \P^{1}$ is smooth and projective. We have $\omega_{X'}=g^{*}\omega_{X}\otimes\O_{X'}(R)$ where $R$ is some effective divisor whose irreducible components are contracted by $g$, hence $g_{*}\omega_{X'}^{m}=\omega_{X}^{m}\otimes g_{*}O(mR)=\omega_{X}^{m}$ for every $m\ge 0$. We may thus replace $X$ with $X'$ and assume that the generic fiber is smooth. This guarantees that $\on{rank}f_{*}\omega_{X}^{m}=\on{rank}f_{*}\omega_{f}^{m}$ has growth $O(m^{\dim X-1})$.
		
		If there are no injective homomorphisms $\O_{\P^{1}}(1)\to f_{*}\omega_{X}^{m}$ for every $m>0$, then $\h^{0}(\omega_{X}^{m})\le\on{rank} f_{*}\omega_{X}^{m}=\on{rank}f_{*}\omega_{f}^{m}$, and this has growth $O(m^{\dim X-1})$.
		
		On the other hand, let $\O_{\P^{1}}(1)\to f_{*}\omega_{X}^{m_{0}}$ be an injective homomorphism for some $m_{0}>0$. In particular, $X$ has Kodaira dimension $\ge 0$.
	
		For some $m$, the closure $Y$ of the image of $X\dashrightarrow \P(\H^{0}(X,\omega_{X}^{mm_{0}}))$ has dimension equal to the Kodaira dimension of $X$ and $k(Y)$ is algebraically closed in $k(X)$, see \cite[\S 3]{iit71}. If $X'$ is a smooth projective variety birational to $X$, then there is a natural isomorphism $\H^{0}(X,\omega_{X}^{mm_{0}})=\H^{0}(X',\omega_{X'}^{mm_{0}})$, see \cite[Ch. 2, Theorem 8.19]{har77}. Thus, up to replacing $X$ with some other smooth, projective variety birational to $X$, we may assume that $X\dashrightarrow Y\s \P(\H^{0}(X,\omega_{X}^{mm_{0}}))$ is defined everywhere and has smooth, projective generic fiber $Z$ by resolution of singularities. Iitaka has then shown that $Z$ has Kodaira dimension $0$, see \cite[Theorem 5]{iit71}. This is easy to see in the case in which $\omega_{X}^{mm_{0}}$ is base point free, since then $\omega_{X}^{mm_{0}}$ is the pullback of $\O(1)$ and thus $\omega_{Z}^{mm_{0}}=\omega_{X}^{mm_{0}}|_{Z}$ is trivial.
		
		Let us recall briefly Grothendieck's convention that, if $V$ is a vector bundle, then $\P(V)$ is the set (or scheme) of linear quotients $V\to k$ up to a scalar. A non-trivial linear map $W\to V$ thus induces a rational map $\P(V)\dashrightarrow\P(W)$ by restriction. If $L$ is a line bundle with non-trivial global sections, the rational map $X\dashrightarrow\P(\H^{0}(X,L))$ is defined by sending a point $x\in X$ outside the base locus to the quotient $\H^{0}(X,L)\to L_{x}\simeq k$. If $L$ embeds in another line bundle $M$, then there is a natural factorization $X\dashrightarrow\P(\H^{0}(X,M))\dashrightarrow \P(\H^{0}(X,L))$, and any point of $X$ outside the support of $M/L$ and outside the base locus of $L$ maps to the locus of definition of $\P(\H^{0}(X,M))\dashrightarrow \P(\H^{0}(X,L))$. 
		
		Let $F\s X$ be the fiber over any rational point of $\P^{1}$. The injective homomorphism $\O_{\P^{1}}(1)\to f_{*}\omega_{X}^{m_{0}}$ induces an injective homomorphism $\O_{\P^{1}}(m)\to f_{*}\omega_{X}^{mm_{0}}$, choose any embedding $\O_{\P^{1}}(1)\to\O_{\P^{1}}(m)$, these induce an injective homomorphism $\O_{X}(F)\to\omega_{X}^{mm_{0}}$. Since $\O_{X}(F)$ induces the morphism $f:X\to\P^{1}$, the composition
		\[X\to Y\s\P(\H^{0}(X,\omega_{X}^{mm_{0}}))\dashrightarrow\P^{1}\]
		coincides with $f$. Observe that the right arrow depends on the choice of the embedding $\O_{X}(F)\to\omega_{X}^{mm_{0}}$, but the composition doesn't.
		
		Let $\xi$ be the generic point of $\P^{1}$, $U\s Y$ an open subset such that $U\to \P^{1}$ is defined, $Y_{\xi}$ the closure of $U_{\xi}$ in $Y$. Then the generic fiber $Z$ of $X\to Y$ is the generic fiber of $X_{\xi}\to Y_{\xi}$, too. By hypothesis, $X_{\xi}$ is of general type, thus by adjunction $\omega_{X_{\xi}}|_{Z}=\omega_{Z}$ is big and hence $Z$ is of general type.
		
		Since $Z$ is a variety of general type of Kodaira dimension $0$ over $\spec k(Y)$, then $Z=\spec k(Y)$, the morphism $X\to Y$ is generically injective and thus $X$ is of general type.
	\end{proof}
\end{proposition}

\begin{remark}
	We don't actually need the precision of \autoref{posgen}: for our purposes it is enough to show that, if $f_{*}\omega_{X}^{m_{0}}$ has a positive \emph{enough} sub-line bundle for some $m_{0}$, then $X$ is of general type. This weaker fact has a more direct proof, let us sketch it. 
	
	First, let us mention an elementary fact about injective sheaf homomorphisms. Let $P,Q$ be vector bundles on $\P^{1}$ and $M,N$ vector bundles on $X$, with $P$ of rank $1$. Suppose we are given injective homomorphisms $m\in\hom(P,f_{*}M)$, $n\in\hom(Q,f_{*}N)$. Then $m^{a}\otimes n\in\hom(P^{\otimes a}\otimes Q,f_{*}(M^{\otimes a}\otimes N))$ is injective for every $a>0$: this can be checked on the generic point of $\P^{1}$ and thus on the generic fiber $X_{k(\P^{1})}$, where the fact that $P$ has rank $1$ allows us to reduce to the fact that the tensor product of non-zero sections of vector bundles is non-zero on an integral scheme.
	
	Assume we have an injective homomorphism $\O_{\P^{1}}(3m_{0})\to f_{*}\omega_{X}^{m_{0}}$, or equivalently $\O_{\P^{1}}(5m_{0})\to f_{*}\omega_{f}^{m_{0}}$, we want to prove that $X$ is of general type. Let $r(m)$ be the rank $f_{*}\omega_{f}^{mm_{0}}$ for every $m$. Since the generic fiber is of general type, up to replacing $m_{0}$ by a multiple $m_{0}'$ we may assume that the growth of $r(m)$ is $O(m^{\dim X-1})$. The induced morphism $\O_{\P^{1}}(5m_{0}')\to f_{*}\omega_{f}^{m_{0}'}$ is injective thanks to the above.
		
	Thanks to \cite[Theorem III]{vie83}, every line bundle in the factorization of $f_{*}\omega_{f}^{mm_{0}}$ has non-negative degree, we may thus choose an injective homomorphism $\O_{\P^{1}}^{r(m)}\to f_{*}\omega_{f}^{mm_{0}}$. Taking the tensor product with the $m$-th power of the homomorphism given by hypothesis, we get an homomorphism $\O_{\P^{1}}(5mm_{0})^{r(m)}\to f_{*}\omega_{f}^{2mm_{0}}$ which is injective thanks to the above.
	
	Since $f_{*}\omega_{X}^{2mm_{0}}=f_{*}\omega_{f}^{2mm_{0}}\otimes\O_{\P^{1}}(-4mm_{0})$, we thus have an injective homomorphism 
	\[\O_{\P^{1}}(mm_{0})^{r(m)}\to f_{*}\omega_{X}^{2mm_{0}}.\]
	In particular, we have $\h^{0}(\omega_{X}^{2mm_{0}})\ge (mm_{0}+1)r(m)$ which has growth $O(m^{n})$, hence $X$ is of general type.
\end{remark}

\begin{corollary}\label{pull}
	Let $f:X\to\P^{1}$ be a non-birationally trivial family of varieties of general type. Then there exists an integer $d_{0}$ and a non-empty open subset $U\s\P^{1}$ such that, for every finite cover $c:\P^{1}\to\P^{1}$ with $\deg c\ge d_{0}$ and such that the branch points of $c$ are contained in $U$, we have that $X_{c}$ is of general type. If $X$ is smooth and projective, $U$ can be chosen as the largest open subset such that $f|_{f^{-1}(U)}$ is smooth.
	\begin{proof}
		By resolution of singularities, we may assume that $X$ is smooth and projective. By generic smoothness, there exists an open subset $U\s\P^{1}$ be such that $f|_{X_{U}}$ is smooth. We have that $X_{c}$ is smooth for every $c:\P^{1}\to\P^{1}$ whose branch points are contained in $U$ since each point of $X_{c}$ is smooth either over $X$ or over $\P^{1}$.
		
		Let $m_{0}$ be the integer given by \autoref{pos}, we have an injective homomorphism $\O(1)\to f_{*}\omega_{f}^{m_{0}}$. Set $d_{0}=2m_{0}+1$, for every finite cover $c$ of degree $\deg c\ge d_{0}=2m_{0}+1$ we have an induced homomorphism $\O(2m_{0}+1)\to f_{c*}\omega_{f_{c}}^{m_{0}}$ and thus $\O(1)\to f_{c*}\omega_{X_{c}}^{m_{0}}$. It follows that $X_{c}$ is of general type thanks to \autoref{posgen}.
	\end{proof}
\end{corollary}

\section{Higher dimensional HIT}

\subsection{Pulling fat sets}

Recall that Serre \cite[Chapter 9]{ser97} defined a subset $S$ of $\P^{1}(k)$ as \emph{thin} if there exists a morphism $f:X\to \P^{1}$ with $X$ of finite type over $k$, finite generic fiber and no generic sections $\spec k(\P^{1})\to X$ such that $S\s f(X(k))$. It's immediate to check that a subset of a thin set is thin, and a finite union of thin sets is thin. Serre's form of Hilbert's irreducibility theorem says that, if $k$ is finitely generated over $\Q$, then $\P^{1}(k)$ is not thin.

\begin{definition}
	A subset $S\s\P^{1}(k)$ is \emph{fat} if the complement $\P^{1}(k)\setminus S$ is thin.
	
	Given a subset $S\s\P^{1}(k)$, a finite set of finite morphisms $D=\{d_{i}:D_{i}\to\P^{1}\}_{i}$ each of degree $>1$ with $D_{i}$ smooth, projective and geometrically connected is a \emph{scale} for $S$ if $S\cup \bigcup_{i}d_{i}(D_{i}(k))=\P^{1}(k)$. The set of branch points of the scale $D$ is the union of the sets of branch points of $d_{i}$.
\end{definition}

Using the fact that a connected scheme with a rational point is geometrically connected \cite[\href{https://stacks.math.columbia.edu/tag/04KV}{Lemma 04KV}]{stacks-project}, it's immediate to check that a subset of $\P^{1}$ is fat if and only if it has a scale. The set of branch points of a scale gives valuable information about a fat set.

\begin{lemma}\label{fat}
	Let $S\s\P^{1}$ be a fat set, and let $D=\{d_{i}:D_{i}\to\P^{1}\}_{i}$ be a scale for $S$. Let $c:\P^{1}\to\P^{1}$ be a morphism such that the sets of branch points of $c$ and $D$ are disjoint. Then $c^{-1}(S)$ is fat.
	\begin{proof}
		Let $d_{i}':D_{i}'\to\P^{1}$ be the base change of $d_{i}$ along $c:\P^{1}\to\P^{1}$. By construction, $c^{-1}(S)\cup\bigcup_{i}d_{i}'(D_{i}'(k))=\P^{1}(k)$. Since the sets of branch points of $c$ and $d_{i}$ are disjoint, we have that $D_{i}'$ is geometrically connected, see for instance \cite[Lemma 2.8]{str20}. Moreover, $D_{i}'$ is smooth since each point of $D_{i}'$ is étale either over $\P^{1}$ or $D_{i}$. It follows that $d_{i}'$ has degree $>1$ and $\{d_{i}':D_{i}'\to\P^{1}\}_{i}$ is a scale for $c^{-1}(S)$, which is thus fat.
	\end{proof}
\end{lemma}

\subsection{Decreasing the fiber dimension}\label{decrease}

Let us now prove Theorem A. Using Hilbert's irreducibility, it's easy to check that Theorem A is equivalent to the following statement.
\vspace{1em}

If the generic fiber of $f:X\to \P^{1}$ is GeM and $f(X(k))$ is fat, there exists a section $\spec k(\P^{1})\to X$.

\vspace{1em}
We prove this statement by induction on the dimension of the generic fiber. If the generic fiber has dimension $0$, this follows from the definition of fat set. Let us prove the inductive step.

We define recursively a sequence of closed subschemes $X_{i+1}\s X_{i}$ with $X_{0}=X$ and such that $f(X_{i}(k))\s\P^{1}_{k}$ is fat.

\begin{itemize}
	\item Define $X'_{i}$ as the closure of $X_{i}(k)$ with the reduced scheme structure, $f(X'_{i}(k))=f(X_{i}(k))\s\P^{1}_{k}$ is fat.
	\item Define $X''_{i}$ as the union of the irreducible components of $X'_{i}$ which dominate $\P^{1}$, $f(X''_{i}(k))\s\P^{1}_{k}$ is fat since $f(X_{i}'(k))\setminus f(X_{i}''(k))$ is finite.
	\item Write $X''_{i}=\bigcup_{j}Y_{i,j}$ as union of irreducible components, $Y_{i,j}\to \P^{1}$ is dominant for every $j$. For every $j$, there exists a finite cover $C_{i,j}\to\P^{1}$ with $C_{i,j}$ smooth projective and a rational map $Y_{i,j}\dashrightarrow C_{i,j}$ with geometrically irreducible generic fiber. If $C_{i,j}\to\P^{1}$ is an isomorphism, define $Z_{i,j}=Y_{i,j}$. Otherwise, there exists a non-empty open subset $V_{i,j}\s Y_{i,j}$ such that $Y_{i,j}\dashrightarrow C_{i,j}$ is defined on $V_{i,j}$. In particular, $f(V_{i,j}(k))\s\P^{1}(k)$ is thin. Define $Z_{i,j}= Y_{i,j}\setminus V_{i,j}$ and $X_{i+1}=\bigcup_{j}Z_{i,j}\s X_{i}$. By construction, $f(X_{i+1}(k))\s\P^{1}(k)$ is fat since $f(X_{i}''(k))\setminus f(X_{i+1}(k))$ is thin.
\end{itemize}

By noetherianity, the sequence is eventually stable, let $r$ be such that $X_{r+1}=X_{r}$. Since $X_{r+1}=X_{r}$, then $X_{r}(k)$ is dense in $X_{r}$, thus every irreducible component is geometrically irreducible, see \cite[\href{https://stacks.math.columbia.edu/tag/0G69}{Lemma 0G69}]{stacks-project}. Moreover, every irreducible component of $X_{r}$ dominates $\P^{1}$ with geometrically irreducible generic fiber. Replace $X$ with $X_{r}$ and write $X=\bigcup_{j}Y_{j}$ as union of irreducible components, we may assume that $Y_{j}\to\P^{1}$ is a family of GeM varieties for every $j$ and $Y_{j}(k)$ is dense in $Y_{j}$.

If $Y_{j}\to\P^{1}$ is birationally trivial for some $j$, since $Y_{j}(k)$ is dense in $Y_{j}$ and a generic fiber of $Y_{j}\to\P^{1}$ has a finite number of rational points, then $\dim Y_{j}=0$, $Y_{j}\to \P^{1}$ is birational and we conclude. Otherwise, thanks to \autoref{pull}, there exists an integer $d_{0}$ and a non-empty open subset $U\s\P^{1}$ such that, for every finite cover $c:\P^{1}\to\P^{1}$ with $\deg c\ge d_{0}$ such that the branch points of $c$ are contained in $U$, we have that $Y_{j,c}$ is of general type for every $j$. 

Let $D=\{d_{l}:D_{l}\to\P^{1}\}$ be a scale for $f(X(k))$. Up to shrinking $U$ furthermore, we may assume that the set of branch points of $D$ is disjoint from $U$. Since we are assuming that the weak Bombieri-Lang conjecture holds up to dimension $\dim X$, the dimension of $\overline{Y_{j,c}(k)}\s Y_{j,c}$ is strictly smaller than $\dim Y_{j}$ for every $j$. Moreover, we have that $f_{c}(X_{c}(k))=m_{c}^{-1}(f(X(k)))$ is fat thanks to \autoref{fat}. It follows that, by induction hypothesis, there exists a generic section $\spec k(\P^{1})\to X_{c}$ for \emph{every} finite cover $c$ as above. There are a lot of such covers: let us show that we can choose them so that the resulting sections "glue" to a generic section $\spec k(\P^{1})\to X$.

\subsection{Gluing sections}\label{glue}

Choose coordinates on $\P^{1}$ so that $0,\infty\in U$, let $p$ be any prime number greater than $d_{0}$. For any positive integer $n$, let $m_{n}:\P^{1}\to\P^{1}$ be the $n$-th power map. We have shown above that there exists a rational section $\P^{1}\dashrightarrow X_{m_{p}}$ for every prime $p\ge d_{0}$, call $s_{p}:\P^{1}\dashrightarrow X_{m_{p}}\to X$ the composition.

We either assume that there exists an integer $N$ such that, for every rational point $v\in \P^{1}(k)$, we have $|X_{v}(k)|\le N$ or that the Bombieri-Lang conjecture holds in every dimension. In the second case, the uniform bound $N$ exists thanks to a theorem of Caporaso-Harris-Mazur and Abramovich-Voloch \cite[Theorem 1.1]{chm97} \cite[Theorem 1.5]{av96} \cite{abr97}. Choose $N+1$ prime numbers $p_{0},\dots, p_{N}$ greater than $d_{0}$, for each one we have a rational section
\[\begin{tikzcd}
														&	X\dar["f"]	\\
	\P^{1}\rar[swap,"m_{p}"]\ar[ur,dashed,"s_{p}"]	&	\P^{1}
\end{tikzcd}\]

Let $Q=\prod_{i=0}^{N}p_{i}$, for every $i=0,\dots,N$, we get a rational section $S_{p_{i}}$ by composition with $s_{p_{i}}$:
\[\begin{tikzcd}[column sep=large]
																									&																		&	X\dar["f"]	\\
	\P^{1}\rar[swap,"m_{Q/p_{i}}"]\ar[urr,dashed,"S_{p_{i}}"]\ar[rr,bend right=35,swap,"m_{Q}"]	&	\P^{1}\rar[swap,"m_{p_{i}}"]\ar[ur,dashed,swap,"s_{p_{i}}"]		&	\P^{1}
\end{tikzcd}\]

Let $V\s\P^{1}$ be an open subset such that $S_{p_{i}}$ is defined on $V$ for every $i$. For every rational point $v\in V(k)$, we have $|X_{v}(k)|\le N$ and thus there exists a couple of different indexes $i\neq j$ such that $S_{p_{i}}(v)=S_{p_{j}}(v)$ for infinitely many $v\in V(k)$, hence $S_{p_{i}}=S_{p_{j}}$. Let $Z\s X$ be the image $S_{p_{i}}=S_{p_{j}}$, by construction we have
\[k(\P^{1})=k(t)\s k(Z)\s k(t^{-p_{i}})\cap k(t^{-p_{j}})\s k(t^{-Q}).\]
Using Galois theory on the cyclic extension $k(t^{-Q})/k(t)$, it is immediate to check that $k(t^{-p_{i}})\cap k(t^{-p_{j}})=k(t)\s k(t^{-Q})$ since $p_{i},p_{j}$ are coprime, thus $k(Z)=k(t)$ and $Z\to\P^{1}$ is birational. This concludes the proof of Theorem A.

\subsection{Non-rational base}

Let us show how Theorem A implies Theorem B. Let $C$ be a geometrically connected curve over a field $k$ finitely generated over $\Q$, and let $f:X\to C$ be a morphism of finite type whose generic fiber is a GeM scheme. Assume that there exists a non-empty open subset $V\s C$ such that $X|_{V}(h)\to V(h)$ is surjective for every finite extension $h/k$. We want to prove that there exists a generic section $C\dashrightarrow X$. It's easy to reduce to the case in which $C$ is smooth and projective, so let us make this assumption.

Observe that, up to replacing $X$ with an affine covering, we may assume that $X$ is affine. Choose $C\to\P^{1}$ any finite map: since $X$ is affine, the Weil restriction $R_{C/\P^{1}}(X)\to\P^{1}$ exists \cite[\S 7.6, Theorem 4]{blr90}. Recall that $R_{C/\P^{1}}(X)\to\P^{1}$ represents the functor on $\P^{1}$-schemes $S\mapsto\hom_{C}(S\times_{\P^{1}}C,X)$.

If $L/k(C)/k(\P^{1})$ is a Galois closure and $\Sigma$ is the set of embeddings $\sigma:k(C)\to L$ as $k(\P^{1})$ extensions, the scheme $R_{C/\P^{1}}(X)_{L}$ is isomorphic to the product $\prod_{\Sigma}X\times_{\spec k(C),\sigma}\spec L$ and hence is a GeM scheme, see \cite[Lemma 3.3]{blfg}. It follows that the generic fiber $R_{C/\P^{1}}(X)_{k(\P^{1})}$ is a GeM scheme, too.

Let $U\s\P^{1}$ be the image of $V\s C$. The fact that $X|_{V}(h)\to V(h)$ is surjective for every finite extension $h/k$ implies that $R_{C/\P^{1}}(X)|_{U}(k)\to U(k)$ is surjective. By Theorem A, we get a generic section $\P^{1}\dashrightarrow R_{C/\P^{1}}(X)$, which in turn induces generic section $C\dashrightarrow X$ by the universal property of $R_{C/\P^{1}}(X)$. This concludes the proof of Theorem B.

\section{Polynomial bijections $\Q\times\Q\to \Q$}

Let us prove Theorem C. Let $k$ be finitely generated over $\Q$, and let $f:\A^{2}\to \A^{1}$ be any morphism. Assume by contradiction that $f$ is bijective on rational points.

First, let us show that the generic fiber of $f$ is geometrically irreducible. This is equivalent to saying that $\spec k(\A^{2})$ is geometrically connected over $\spec k(\A^{1})$, or that $k(\A^{1})$ is algebraically closed in $k(\A^{2})$. Let $k(\A^{1})\s L\s k(\A^{2})$ a subextension algebraic over $k(\A^{1})$. Let $C\to\A^{1}$ be a finite cover with $C$ regular and $k(C)=L$. The rational map $\A^{2}\dashrightarrow C$ is defined in codimension $1$, thus there exists a finite subset $S\s \A^{2}$ and an extension $\A^{2}\setminus S\to C$. Since the composition $\A^{2}\setminus S(k)\to C(k)\to\A^{1}(k)$ is surjective up to a finite number of points, by Hilbert's irreducibility theorem we have that $C=\A^{1}$, i.e. $L=k(\A^{1})$.

This leaves us with three cases: the generic fiber is a geometrically irreducible curve of geometric genus $0$, $1$, or $\ge 2$. The first two have been settled by W. Sawin in the polymath project \cite{tao19}, while the third follows from Theorem A. Let us give details for all of them.

\subsection{Genus 0} Assume that the generic fiber of $f$ has genus $0$. By generic smoothness, there exists an open subset $U\s \A^{2}$ such that $f|_{U}$ is smooth. For a generic rational point $u\in U(k)$, the fiber $f^{-1}(f(u))$ is birational to a Brauer-Severi variety of dimension $1$ and has a smooth rational point, thus it is birational to $\P^{1}$ and $f^{-1}(f(u))(k)$ is infinite. This is absurd.

\subsection{Genus 1} Assume now that the generic fiber has genus $1$. By resolution of singularities, there exists an open subset $V\s\A^{1}$, a variety $X$ with a smooth projective morphism $g:X\to V$ whose fibers are smooth genus $1$ curves and a compatible birational map $X\dashrightarrow\A^{2}$. Up to shrinking $V$, we may suppose that the fibers of $f|_{V}$ are geometrically irreducible. Let $U$ be a variety with open embeddings $U\s X$, $U\s\A^{2}$, replace $V$ with $g(U)\s V$ so that $g|_{U}$ is surjective. 

The morphism $X\setminus U\to V$ is finite, let $N$ be its degree. Since the fibers of $U\to V$ have at most one rational point, it follows that $|X_{v}(k)|\le N+1$ for every $v\in V(k)$.

Every smooth genus $1$ fibration is a torsor for a relative elliptic curve (namely, its relative $\underline{\Pic}^{0}$), thus there exists an elliptic curve $E\to V$ such that $X$ is an $E$-torsor. Moreover, every torsor for an abelian variety is torsion, thus there exists a finite morphism $\pi:X\to E$ over $V$ induced by the $n$-multiplication map $E\to E$ for some $n$.

If $v\in V(k)$ is such that $X_{v}(k)$ is non-empty, then $|X_{v}(k)|=|E_{v}(k)|\le N+1$. This means that, up to composing $\pi$ with the $(N+1)!$ multiplication $E\to E$, we may assume that $\pi(X(k))\s V(k)\s E(k)$, where $V\to E$ is the identity section. In particular, $X(k)\s\pi^{-1}(V(k))$ is not dense. This is absurd, since $X$ is birational to $\A^{2}$.

\subsection{Genus $\ge 2$}

Thanks to Theorem A, there exists an open subset $V\s\A^{1}$ and a section $s:V\to\A^{2}$. It follows that $\A^{2}|_{V}(k)=s(V(k))$, which is absurd since $s(V)$ is a proper closed subset and $\A^{2}|_{V}(k)$ is dense.

\printbibliography
	
\end{document}